\documentclass[oneside]{amsart}

\usepackage{graphicx}%
\usepackage{multirow}%
\usepackage{amsmath,amssymb,amsfonts}%
\usepackage{amsthm}%
\usepackage{tikz}

\theoremstyle{plain}
\newtheorem{thm}{Theorem}%[section]
\newtheorem{cor}[thm]{Corollary}
\newtheorem{defn}[thm]{Definition}
\newtheorem{prop}[thm]{Proposition}

\DeclareMathOperator{\ex}{exc}
\DeclareMathOperator{\Ex}{Exc}
\DeclareMathOperator{\df}{dfc}

\DeclareMathOperator{\diag}{diag}
\DeclareMathOperator{\adiag}{adiag}

\begin{document}

\title{Excedances and Deficiencies of Integer Partitions}

\author{George Andrews and Brian Hopkins}

\begin{abstract}
The notion of a fixed point, applied by Blecher and Knopfmacher to integer partitions, can be considered as a permutation statistic.  Here we apply two related permutation statistics, excedance and deficiency, to integer partitions.  Using analytic and combinatorial techniques, we establish several generating functions and make connections to partition concepts including Dyson's crank statistic, overpartitions, and the Rogers--Ramanujan partitions.
\end{abstract}


\maketitle

\section{Introduction and background}

A partition of a positive integer $n$ is a weakly decreasing sequence of positive integers whose sum
is $n$. The positive integers in the sequence are called parts.  Let $P(n)$ be the set of partitions of $n$ and write $p(n) = | P(n) |$.  By convention, $p(0) = 1$.

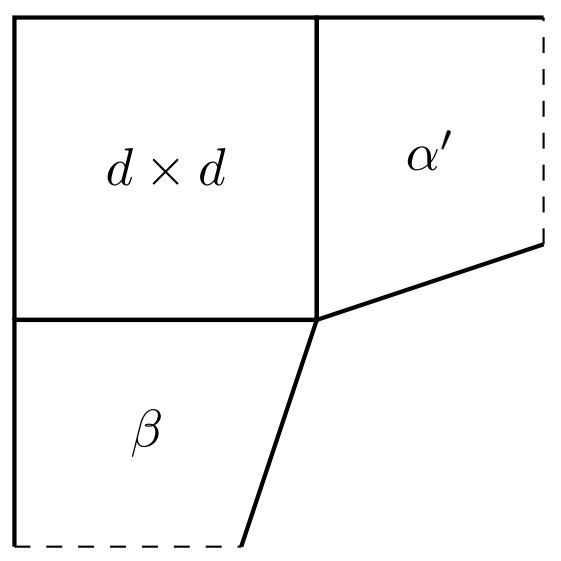
\begin{figure}[th]
 \centering
\begin{tikzpicture}[scale=0.5]
\draw[thick] (0,0) -- (0,3) -- (4,3) -- (3,0);
\draw[thick] (0,3) -- (0,7) -- (4,7) -- (4,3) -- (0,3);
\draw[thick] (7,4) -- (4,3) -- (4,7) -- (7,7);
\draw[dashed] (7,4) -- (7,7);
\draw[dashed] (0,0) -- (3,0);
\node at (2,5) {$d \times d$};
\node at (1.75,1.5) {$\beta$};
\node at (5.5,5.25) {$\alpha'$};
\end{tikzpicture}
 \caption{The Durfee square decomposition of a partition.} \label{dab}
 \end{figure}

The Durfee square of a partition $\lambda = (\lambda_1, \ldots, \lambda_t)$ has dimension $d$ where $d$ is the greatest index for which $\lambda_d \ge d$.  This statistic gives a combinatorial understanding of Jacobi's result
\begin{equation}
\sum_{n \ge 0} p(n)q^n =  \frac{1}{(q;q)_\infty} = \sum_{d \ge 0} \frac{q^{d^2}}{(q;q)_d^2} \label{jac}
 \end{equation}
 where we use standard $q$-series notation; see, for example, Chapter 10 of Andrews, Askey, and Roy \cite{aar}.  As shown in Figure \ref{dab}, a partition can be decomposed into its $d \times d$ Durfee square, the conjugate $\alpha'$ of a partition $\alpha$ with largest part at most $d$ to the right of the Durfee square, and a partition $\beta$ with largest part at most $d$ below the Durfee square.
 
 We will also use the ${}_2 \phi_1$ series studied by Heine, in particular the identity \cite[Equation 10.10.1]{aar}:
 \begin{equation}
 {}_2 \phi_1 \!\binom{a,b;q,x}{c} = \frac{(c/b; q)_\infty (bx; q)_\infty}{(c;q)_\infty (x;q)_\infty}\, {}_2 \phi_1 \!\binom{abx/c, b; q, c/b}{bx} \label{10101}.
 \end{equation}

The crank of a partition, named by Dyson in 1944 but not found until 1988 work of Andrews and Garvan \cite{ag}, is arguably the most important partition statistic.  We will not need its definition, but will use the following generating function for the number of partitions of $n$ with positive crank.  This was established in both Andrews--Newman \cite[Theorem 1]{an20} and Hopkins--Sellers \cite[proof of Theorem 1]{hs20},  2020 works that independently established a connection between the crank and another partition statistic called the mex.

\begin{prop} \label{nonnegcrankgf}
A generating function for partitions with non-negative crank is
\[ \frac{1}{(q;q)_\infty} \sum_{n \ge 0} (-1)^n q^{\binom{n+1}{2}}.\]
\end{prop}

In 2022, Blecher and Knopfmacher applied the pervasive mathematical notion of fixed points to partitions as follows \cite{bk22}.

\begin{defn}
A partition $\lambda = (\lambda_1, \ldots, \lambda_t)$ has a fixed point if $\lambda_i = i$ for some index $i$.
\end{defn}

Note that, due to the ordering convention, a partition has at most one fixed point.  Also, if $\lambda_i = i$, then $i$ is the size of the Durfee square.

In resolving a conjecture of Blecher and Knopfmacher, Hopkins and Sellers established the following connections between fixed points and the crank \cite{hs24}.

\begin{prop} \label{fixedposcrank}
Given $n \ge 1$, the number of partitions of $n$ with a fixed point equals the number of partitions of $n$ with positive crank.
\end{prop}

To refine the notion of fixed point, let $f(n,d)$ be the number of partitions $\lambda$ of $n$ that have a fixed point $d$, i.e., $\lambda_d  = d$.  Hopkins established the following recurrence \cite{h24}.

\begin{prop} \label{fnd}
For $d \ge 2$ and $n \ge d^2$,
\[ f(n,d) = f(n-1+1,d) + f(n-d,d) - f(n-2d+1,d) + f(n-2d+1,d-1).\]
\end{prop}

Considering fixed points in the context of permutations, partitions without a fixed point correspond to derangements.  This setting invites the application of other permutation statistics to integer partitions.  In this work, we consider two related permutation statistics on partitions.

\begin{defn}
An excedance of a partition $\lambda \in P(n)$ is an index $i$ such that $\lambda_i > i$ while a deficiency is an index $i$ such that $\lambda_i < i$.  For $k \ge 0$, let $\ex(n,k)$ be the number of partitions of $n$ with $k$ excedances, similarly $\df(n,k)$ for deficiencies.
\end{defn}

The terminology dates at least to Reilly 1979 \cite{r}.  Note that every part of a partition is an excedance, fixed point, or deficiency.  Further, any excedances are the initial parts of a partition, followed by a possible fixed point, with any deficiencies comprising the subsequent parts.

The next two sections give results on excedances and deficiencies, respectively, using both generating function and combinatorial techniques.

\section{Excedances of integer partitions}
The array of Figure \ref{exctri} gives the number of partitions of $n$ with $k$ excedances for small values of $n$ and $k$.

\begin{figure}[b]
\begin{tabular}{r|rrrrr}
$n \backslash k$ & 0 & 1 & 2 & 3 & 4 \\ \hline
0 & 1 \\
1 & 1 & & &  \\
2 & 1 & 1 & &  \\
3 & 1 & 2 & &   \\
4 & 1 & 4 & &  \\
5 & 1 & 6 & & \\
6 & 1 & 9 & 1 &  \\
7 & 1 & 12 & 2 &  \\
8 & 1 & 16 & 5 &  \\
9 & 1 & 20 & 9 &  \\
10 & 1 & 25 & 16 &  \\
11 & 1 & 30 & 25 &  \\
12 & 1 & 36 & 39 & 1  \\
13 & 1 & 42 & 56 & 2 \\
14 & 1 & 49 & 80 & 5 \\
15 & 1 & 56 & 109 & 10 \\
16 & 1 & 64 & 147 & 19 \\
17 & 1 & 72 & 192 & 32 \\
18 & 1 & 81 & 249 & 54 \\
19 & 1 & 90 & 315 & 84 \\
20 & 1 & 100 & 396 & 129 & 1
\end{tabular}
\caption{Array of $\ex(n,k)$ values for $0 \le n \le 20$ and $0 \le k \le 4$.} \label{exctri}
\end{figure}

\begin{thm} \label{excthm}
A bivariate generating function for the number of partitions of $n$ with $k$ excedances is 
\begin{align}
 \sum_{n,k} \ex(n,k)w^kq^n & = 1+ \sum_{d \ge 1} \frac{w^{d-1}q^{d^2}(1-q^d+wq^d)}{(q;q)_d^2}. \label{exc1} \\
 & = \frac{1}{(q;q)_\infty} \sum_{n \ge 0} \frac{(wq;q)_n (-1)^n q^{\binom{n+1}{2}}(1+q^{n+1})}{(q;q)_n}. \label{exc2}
 \end{align}
\end{thm}

\begin{proof}
A partition $\lambda$ with Durfee square size $d$ has $d-1$ excedances if $\lambda_d = d$, that is, referencing Figure \ref{dab}, if $d$ is not a part of the subpartition $\alpha'$ to the right of the Durfee square.  Otherwise, $\lambda$ has $d$ excedances and $d$ is a part of $\alpha'$ (this case includes the empty partition of 0).  A generating function for $\ex(n,k)$ is then
\begin{equation}
 \sum_{d \ge 1} \frac{w^{d-1}q^{d^2}}{(q;q)_{d-1}(q;q)_d} + \sum_{d \ge 0} \frac{w^d q^{d^2+d}}{(q;q)_d^2}  \label{exc0}
 \end{equation}
which simplifies to the desired expression \eqref{exc1}.

Working from each term of \eqref{exc0},
\begin{align*}
\sum_{d \ge 1} \frac{w^{d-1}q^{d^2}}{(q;q)_{d-1}(q;q)_d} & = \sum_{d \ge 0} \frac{w^d q^{d^2+2d+1}}{(q;q)_d (q;q)_{d+1}}  \\
& = \frac{q}{1-q} \lim_{\tau \to 0} {}_2\phi_1 \! \binom{q/\tau, q/\tau; q, \tau^2 q w}{q^2} \\
& = \frac{q}{1-q} \frac{1}{(q^2;q)_\infty} \lim_{\tau \to 0} {}_2\phi_1 \! \binom{wq, q/\tau; q, \tau q}{\tau q^2 w} \\
& = \frac{1}{(q;q)_\infty} \sum_{n \ge 0} \frac{(wq;q)_n (-1)^n q^{\binom{n+2}{2}}}{(q;q)_n},
\end{align*}
where we used \eqref{10101} and, similarly,
\begin{align*}
\sum_{d \ge 0} \frac{w^d q^{d^2+d}}{(q;q)^2_d} & = \lim_{\tau \to 0} {}_2\phi_1 \! \binom{q/\tau, q/\tau; q, \tau^2 w}{q} \\
& = \frac{1}{(q;q)_\infty} \lim_{\tau \to 0} {}_2\phi_1 \! \binom{wq, q/\tau; q, \tau}{\tau q w} \\
& = \frac{1}{(q;q)_\infty} \sum_{n \ge 0} \frac{(wq;q)_n (-1)^n q^{\binom{n+1}{2}}}{(q;q)_n},
\end{align*}
which together give \eqref{exc2}.
\end{proof}

The results of Theorem \ref{excthm} allow us to explain many patterns in the $\ex(n,k)$ array of Table \ref{exctri}.

Clearly the row sums of the $\ex(n,k)$ array match $p(n)$, the number of partitions of $n$.  Analytically, letting $w = 1$ in \eqref{exc1} gives Jacobi's \eqref{jac}.

\begin{cor}
A generating function for the number of partitions of $n$ with an even number of excedances is
\[ \frac{1}{2(q;q)_\infty} \sum_{n \ge 0} \left(1 - \frac{(-q;q)_n}{(q;q)_n}\right) (-1)^n q^{\binom{n+1}{2}} (1+q^{n+1})\]
and an analogue for odd number of excedances is
\[ \frac{1}{2(q;q)_\infty} \sum_{n \ge 0} \left(1 + \frac{(-q;q)_n}{(q;q)_n}\right) (-1)^n q^{\binom{n+1}{2}} (1+q^{n+1})\]
giving an expression for the difference
\[ \frac{1}{(q;q)_\infty} \sum_{n \ge 0} \frac{(-q;q)_n (-1)^n q^{\binom{n+1}{2}}(1+q^{n+1})}{(q;q)_n}. \]
\end{cor}
\begin{proof}  
The excedance parity results follow from standard techniques applied to \eqref{exc2} and setting $w = 1$.  The difference expression comes most easily from setting $w = -1$ in \eqref{exc2}.
\end{proof}

Table \ref{excpar} provides initial values for these counts.  Computations suggest that the difference remains positive for $n \ge 12$.  At the time of writing, these sequences were not contained in the On-Line Encyclopedia of Integer Sequences \cite{o}.

\begin{table}[ht]
\caption{The number of partitions of $n$ with an even number of excedances, an odd number, and the difference, through $n = 15$.} \label{excpar}
\begin{tabular}{r|rrrrrrrrrrrrrrr}
$n$ & 1 & 2 & 3 & 4 & 5 & 6 & 7 & 8 & 9 & 10 & 11 & 12 & 13 & 14 & 15 \\ \hline
even & 1 & 1 & 1 & 1 & 1 & 2 & 3 & 6 & 10 & 17 & 26 & 40 & 57 & 81 & 110 \\
odd & 0 & 1 & 2 & 4 & 6 & 9 & 12 & 16 & 20 & 25 & 30 & 37 & 44 & 54 & 66 \\
diff. & 1 & 0 & $-1$ & $-3$ & $-5$ & $-7$ & $-9$ & $-10$ & $-10$ & $-8$ & $-4$ & 3 & 13 & 27 & 44 
\end{tabular}
\end{table}

The weighted row sum $\ex(n) = \sum_k k \ex(n,k)$ gives the total number of excedances in the partitions of $n$.

\begin{cor}
A generating function for the total number of excedances in the partitions of $n$ is
\[ \frac{1}{(q;q)_\infty} \sum_{n \ge 1} \frac{(-1)^{n-1} q^{n(n+3)/2}}{1-q^n}.\]
\end{cor}
\begin{proof}
The result follows from determining the partial derivative of \eqref{exc2} with respect to $w$ and evaluating that expression at $w = 1$.
\end{proof}

Values for $\ex(n)$ for small $n$ are given in Table \ref{toted} in the next section.

Both the diagonal sums 
\[\diag(\ex,n) = \ex(n,0) + \ex(n-1,1) + \cdots \]
and antidiagonal sums
\[\adiag(\ex,n) = \ex(n,0) + \ex(n+1, 1) + \cdots \]
(finite due to the shape of the array) have combinatorial interpretations.

\begin{cor} \label{excdiag}
For $n \ge 1$, the diagonal sum $\diag(\ex,n)$ equals the number of partitions of $n+1$ with positive crank.
\end{cor}

\begin{proof}
Letting $w = q$ in \eqref{exc2} gives the diagonal sums of the $\ex(n,k)$ array.  The resulting expression reduces to
\[ \frac{1}{(q;q)_\infty} \sum_{n \ge 1} (-1)^{n+1} q^{\binom{n+1}{2}}\]
which is the complement of the shifted generating function of Proposition \ref{nonnegcrankgf}, i.e., the expression here is a generating function for the number of partitions of $n+1$ with positive crank, using crank symmetry \cite{ag}.
\end{proof}

\begin{cor} \label{excadiag}
The antidiagonal sum $\adiag(\ex,n)$ equals $p(n) + p(n-1)$.
\end{cor}

\begin{proof}
Letting $w = 1/q$ in \eqref{exc2} gives the antidiagonal sums of the $\ex(n,k)$ array.  The resulting term $(1;q)_n$ in the numerator is nonzero only for $n=0$ so that the expression reduces to
\[\frac{1+q}{(q;q)_\infty},\]
a generating function for $p(n) + p(n-1)$.
\end{proof}

We also give combinatorial proofs for Corollaries \ref{excdiag} and \ref{excadiag}.  Let $\Ex(n,k)$ be the set of partitions of $n$ with $k$ excedances.
Write $F(n,d)$ for the set of partitions $\lambda \in P(n)$ with a fixed point $\lambda_d = d$ and $F(n)$ for the set of partitions in $P(n)$ with a fixed point \cite{h24}.

\begin{thm} \label{ef}
There is a bijection
\[\Ex(n,k) \cong F(n+k+1,k+1).\]
\end{thm}
\begin{proof}
Given $\lambda \in \Ex(n,k)$, the $k$ excedances are necessarily the first $k$ parts.  Note that $\lambda_k \ge k+1$ and $\lambda_{k+1} \le k+1$ (else there are at least $k+1$ excedances).  
Assign 
\[\lambda \mapsto \mu = (\lambda_1, \ldots, \lambda_k, k+1, \lambda_{k+1}, \ldots ),\]
a partition of $n+k+1$ with fixed point $\mu_{k+1} = k+1$ so that $\mu \in F(n+k+1,k+1)$.  

The reverse map is clear:\ remove the fixed point from $\mu \in F(n+k+1,k+1)$ to get $\lambda \in \Ex(n,k)$.
\end{proof}

Examples of the bijection of Theorem \ref{ef} are given in Table \ref{9ex}.

\begin{proof}[Combinatorial proof of Corollary \ref{excdiag}]  By Theorem \ref{ef}, the diagonal sum
\begin{align*} \diag(\ex,n) & = \sum_{i \ge 0} \ex(n-i,i) \\
& = \sum_{i \ge 0} f(n+1,i+1) \\
& = \sum_{j \ge 1} f(n+1,j),
\end{align*}
the number of partitions of $n+1$ with a fixed point.  By Proposition \ref{fixedposcrank}, this matches the number of partitions of $n+1$ with positive crank.
\end{proof}

Table \ref{9ex} gives an example of the bijection in this combinatorial proof of Corollary \ref{excdiag} which uses three applications of the bijection of Theorem \ref{ef}.

\begin{table}[b]
\caption{The $n=8$ example of the bijection in the second proof of Corollary \ref{excadiag}.} \label{9ex}
\begin{tabular}{rr|rr}
& $\diag(\Ex,8)$ & $F(9)$ &  \\ \hline \hline
$\Ex(8,0)$ & $(1,1,1,1,1,1,1,1)$ & $(1,1,1,1,1,1,1,1,1)$ & $F(9,1) $\\ \hline
$\Ex(7,1)$ & $(7)$ & $(7,2)$ & $F(9,2)$ \\
& $(6,1)$ & $(6,2,1)$ & \\
& $(5,2)$ & $(5,2,2)$ & \\
& $(5,1,1)$ & $(5,2,1,1)$ & \\
& $(4,2,1)$ & $(4,2,2,1)$ & \\
& $(4,1,1,1)$ & $(4,2,1,1,1)$ & \\
& $(3,2,2)$ & $(3,2,2,2)$ & \\
& $(3,2,1,1)$ & $(3,2,2,1,1)$ & \\
& $(3,1,1,1,1)$ & $(3,2,1,1,1,1)$ & \\
& $(2,2,2,1)$ & $(2,2,2,2,1)$ & \\
& $(2,2,1,1,1)$ & $(2,2,2,1,1,1)$ \\
& $(2,1,1,1,1,1)$ & $(2,2,1,1,1,1)$ & \\ \hline
$\Ex(6,2)$ & $(3,3)$ & $(3,3,3)$ & $F(9,3)$ \\ \hline
\end{tabular}
\end{table}

Theorem \ref{ef} also translates Proposition \ref{fnd} into a recurrence for the $\ex(n,k)$ array.

\begin{cor}
For $k \ge 1$ and $n \ge k^2+k$,
\[ \ex(n,k) = \ex(n-k,k) + \ex(n-k-1,k) - \ex(n-2k-1,k) + \ex(n-2k,k-1).\]
\end{cor}

For the next combinatorial proof, note that $\lambda \in \Ex(n,d)$ either has Durfee square size $d$ with no fixed point, i.e., $\lambda_d \ge d+1$,  or Durfee square size $d+1$ with fixed point $\lambda_{d+1} = d+1$.

\begin{proof}[Combinatorial proof of Corollary \ref{excadiag}]  We establish a bijection
\[\bigcup_{i \ge 0} \Ex(n+i,k+i) \cong P(n) \cup P(n-1).\]

For $\lambda \in \Ex(n,d)$ with Durfee square size $d$, we know $\lambda_d \ge d+1$ and $\lambda_{d+1} \le d$ so that reducing each excedance by one gives a valid partition
\[ (\lambda_1 - 1, \ldots, \lambda_d - 1, \lambda_{d+1}, \ldots ) \in P(n-d).\]

For $\lambda \in \Ex(n,d)$ with Durfee square size $d+1$, we know $\lambda_{d+1} = d+1$ and removing that fixed point leaves
\[ (\lambda_1, \ldots, \lambda_d, \lambda_{d+2}, \ldots \in P(n-d-1).\]

Applying the map to $\bigcup \Ex(n+i,k+i)$, the first case gives partitions in $P(n)$ and the second case gives partitions in $P(n-1)$.

For the reverse map, given $\mu \in P(n)$ with Durfee square size $d$, increase the first $d$ parts by one to make $\lambda \in P(n+d)$ with $d$ excedances, Durfee square size $d$, and no fixed point.  Given $\nu \in P(n-1)$, find the greatest $e$ such that $\nu_{e-1} \ge e$ and add a part $e$ to make a partition in $P(n+e-1)$ with fixed point $e$ and $e-1$ excedances.  The bijection establishes the identity $\adiag(\ex,n) = p(n) + p(n-1)$.
\end{proof}

Table \ref{10ex} gives an example of the bijection in this combinatorial proof of Corollary \ref{excadiag}.

Note that, by adding an overlined part 1 to the partitions of $P(n-1)$, Corollary \ref{excadiag} can be cast as a bijection to the restricted overpartitions of $n$ \cite{cl} where only a part 1 may be overlined.

\begin{table}[b]
\caption{The $n=5$ example of the bijection in the second proof of Corollary \ref{excadiag}.} \label{10ex}
\begin{tabular}{rr|rr}
& $\adiag(\Ex,5)$ & $P(5)$ & ${} \cup P(4)$ \\ \hline \hline
$\Ex(5,0)$ & $(1,1,1,1,1)$ & & $(1,1,1,1)$ \\ \cline{1-2}
$\Ex(6,1)$ & $(6)$ & $(5)$ & \\
& $(5,1)$ & $(4,1)$ & \\
& $(4,2)$ & & $(4)$ \\
& $(4,1,1)$ & $(3,1,1)$ & \\
& $(3,2,1)$ & & $(3,1)$ \\
& $(3,1,1,1)$ & $(2,1,1,1)$ & \\
& $(2,2,2)$ & & $(2,2)$ \\
& $(2,2,1,1)$ & & $(2,1,1)$ \\
& $(2,1,1,1)$ & $(1,1,1,1,1)$ & \\ \cline{1-2}
$\Ex(7,2)$ & $(4,3)$ & $(3,2)$ & \\
& $(3,3,1)$ & $(2,2,1)$ & \\ \hline
\end{tabular}
\end{table}

\section{Deficiencies of integer partitions}

For permutations, excedances and deficiencies are essentially equivalent, since each excedance of a permutation $\pi$ is a deficiency of the inverse permutation $\pi^{-1}$ and vice versa.  That means that the number of excedances over all permutation of $S_n$ equals the number of deficiencies over $S_n$.  The situation for partitions is very different.  As explained in the proof of Theorem \ref{excthm}, a partition of $n$ with Durfee square size $d$ has $d-1$ or $d$ excedances while the number of deficiencies, the number of parts in $\beta$ below the Durfee square as in Figure \ref{dab}, could be as high as $n - d^2$.  See Corollary \ref{dfctot} and Table \ref{toted} for more details.

The array of Figure \ref{dfctri} gives the number of partitions of $n$ with $k$ deficiencies for small values of $n$ and $k$. 

\begin{figure}[b]
\begin{tabular}{r|rrrrrrrrrrrrrrr}
$n \backslash k$ & 0 & 1 & 2 & 3 & 4 & 5 & 6 & 7 & 8 & 9 & 10 & 11 & 12 & 13 & 14  \\ \hline
1 & 1 & & &  \\
2 & 1 & 1 & &  \\
3 & 1 & 1 & 1 &   \\
4 & 2 & 1 & 1 &  1 \\
5 & 2 & 2 & 1 & 1 & 1 \\
6 & 3 & 3 & 2 & 1 & 1& 1 \\
7 & 3 & 4 & 3 & 2 & 1 & 1 & 1  \\
8 & 4 & 5 & 5 & 3 & 2 & 1 & 1 & 1  \\
9 & 5 & 6 & 6 & 5 & 3 & 2 & 1 & 1 & 1  \\
10 & 6 & 8 & 8 & 7 & 5 & 3 & 2 & 1 & 1 & 1 \\
11 & 7 & 10 & 10 & 9 & 7 & 5 & 3 & 2 & 1 & 1 & 1 \\
12 & 9 & 13 & 13 & 12 & 10 & 7 & 5 & 3 & 2 & 1 & 1 & 1 \\
13 &  10 & 16 & 17 & 15 & 13 & 10 & 7 & 5 & 3 & 2 & 1 & 1 & 1 \\
14 & 12 & 20 & 22 & 20 & 17 & 14 & 10 & 7 & 5 & 3 & 2 & 1 & 1 & 1 \\
15 & 14 & 24 & 28 & 26 & 22 & 18 & 14 & 10 & 7 & 5 & 3 & 2 & 1 & 1 & 1
\end{tabular}
\caption{Array of $\df(n,k)$ values for $0 \le n, k \le 15$.} \label{dfctri}
\end{figure}

\begin{thm} \label{dfcthm}
A bivariate generating function for the number of partitions of $n$ with $k$ deficiencies is 
\begin{align}
 \sum_{n,k} \df(n,k)w^kq^n & = \sum_{d \ge 0} \frac{q^{d^2}}{(q;q)_d (wq;q)_d} \label{dfc1} \\
 & = \frac{1}{(wq;q)_\infty} \sum_{n \ge 0} \frac{(w^{-1};q)_n (-1)^n w^n q^{\binom{n+1}{2}}}{(q;q)_n}. \label{dfc2}
 \end{align}
\end{thm}

\begin{proof}
The number of deficiencies of a partition $\lambda$ is exactly the number of parts of the subpartition $\beta$ underneath the Durfee square, which establishes \eqref{dfc1}.

Similar to the derivations of Theorem \ref{excthm}, we use \eqref{10101} here:
\begin{align*}
\sum_{d \ge 0} \frac{q^{d^2}}{(q;q)_d (wq;q)_d} & = \lim_{\tau \to 0} {}_2\phi_1 \! \binom{q/\tau,1/\tau; q, \tau^2}{wq} \\
& = \lim_{\tau \to 0} \frac{(w \tau;q)_\infty (\tau;q)_\infty}{(wq;q)_\infty (\tau^2;q)_\infty} \, {}_2\phi_1 \! \binom{1/w, 1/\tau; q, wq \tau}{\tau} \\
& = \frac{1}{(wq;q)_\infty} \sum_{n \ge 0} \frac{(w^{-1};q)_n (-1)^n w^n q^{\binom{n+1}{2}}}{(q;q)_n},
\end{align*}
as desired.
\end{proof}

It is interesting to contrast the rather complicated expressions \eqref{exc2} and \eqref{dfc2} to
\[ \sum_{d \ge 0} \frac{w^d q^{d^2}}{(q;q)_d (wq;q)_d} = \frac{1}{(wq;q)_\infty} \]
where each part is weighted by $w$, not just the excedances or deficiencies.

Similar to the previous section, the results of Theorem \ref{dfcthm} allow us to explain many patterns in the $\df(n,k)$ array of Table \ref{dfctri}.

\begin{cor}
The sequence $\df(n,0)$ equals the number of Rogers--Ramanujan partitions of $n$, those restricted to parts congruent to 1 or 4 modulo 5.
\end{cor}

\begin{proof}
Removing the $(wq;q)_d$ term from the denominator in \eqref{dfc1} leaves a generating function for the $w=0$ sequence, and the resulting 
\[\sum_{d \ge 0} \frac{q^{d^2}}{(q;q)_d}\]
is a known generating function for the described Rogers--Ramanujan partitions \cite[p. 52]{ae}.
\end{proof} 

\begin{cor} \label{cor:dfcpar}
A generating function for the number of partitions of $n$ with an even number of deficiencies is
\[ \frac{1}{2(q;q)_\infty}  - \frac{1}{(-q;q)_\infty} \sum_{n \ge 1} \frac{(-q;q)_{n-1} q^{\binom{n+1}{2}}}{(q;q)_n}\]
and an analogue for odd number of deficiencies is
\[ \frac{1}{2(q;q)_\infty}  + \frac{1}{(-q;q)_\infty} \sum_{n \ge 1} \frac{(-q;q)_{n-1} q^{\binom{n+1}{2}}}{(q;q)_n}\]
giving an expression for the difference
\[ \frac{2}{(-q;q)_\infty} \sum_{n \ge 1} \frac{(-q;q)_{n-1} q^{\binom{n+1}{2}}}{(q;q)_n}. \]
\end{cor}
\begin{proof}  
The deficiencies parity results follow from standard techniques applied to \eqref{dfc2} and setting $w = 1$.  The difference expression comes most easily from setting $w = -1$ in \eqref{dfc2}.
\end{proof}

Table \ref{dfcpar} provides initial values for the sequences described in Corollary \ref{cor:dfcpar}. The even counts match \cite[A046682]{o} and odd counts match \cite[A000701]{o}, while the difference matches the number of self-conjugate partitions \cite[A000700]{o}.

\begin{table}[ht]
\caption{The number of partitions of $n$ with an even number of deficiencies, an odd number, and the difference, through $n = 15$.} \label{dfcpar}
\begin{tabular}{r|rrrrrrrrrrrrrrr}
$n$ & 1 & 2 & 3 & 4 & 5 & 6 & 7 & 8 & 9 & 10 & 11 & 12 & 13 & 14 & 15 \\ \hline
even & 1 & 1 & 2 & 3 & 4 & 6 & 8 & 12 & 16 & 22 & 29 & 40 & 52 & 69 & 90 \\
odd & 0 & 1 & 1 & 2 & 3 & 5 & 7 & 10 & 14 & 20 & 27 & 37 & 49 & 66 & 86 \\
diff. & 1 & 0 & 1 & 1 & 1 & 1 & 1 & 2 & 2 & 2 & 2 & 3 & 3 & 3 & 4
\end{tabular}
\end{table}

The weighted row sum $\df(n) = \sum_k k \df(n,k)$ gives the total number of deficiencies in the partitions of $n$

\begin{cor} \label{dfctot}
A generating function for the total number of deficiencies in the partitions of $n$ is
\[ \frac{1}{(q;q)_\infty} \sum_{n \ge 1} \frac{q^n+(-1)^{n} q^{\binom{n+1}{2}}}{1-q^n}.\]
\end{cor}
\begin{proof}
The result follows from determining the partial derivative of \eqref{dfc2} with respect to $w$ and evaluating that expression at $w = 1$.
\end{proof}

Table \ref{toted} gives values of $\ex(n)$ and $\df(n)$ for small $n$, informing the discussion at the beginning of this section.  The $\df(n)$ sequence is \cite[A114089]{o} (with this interpretation); the $\ex(n)$ sequence is not in the OEIS at the time of writing.

\begin{table}[th]
\caption{The total number of excedances and deficiencies in $P(n)$ for small $n$.} \label{toted}
\begin{tabular}{r|rrrrrrrrrrrrrrr}
$n$ & 1 & 2 & 3 & 4 & 5 & 6 & 7 & 8 & 9 & 10 & 11 & 12 & 13 & 14 & 15 \\ \hline
$\ex(n)$ & 0 & 1 & 2 & 4 & 6 & 11 & 16 & 26 & 38 & 57 & 80 & 117 & 160 & 224 & 304 \\
$\df(n)$ & 0 & 1 & 3 & 6 & 11 & 19 & 31 & 50 & 76 & 116 & 169 & 247 & 349 & 494 & 682
\end{tabular}
\end{table}

Again, the diagonal sums 
\[\diag(\df,n) = \df(n,0) + \df(n-1,1) + \cdots \]
have a combinatorial interpretation (here the antidiagonal sums are infinite).

\begin{cor} \label{dfcdiag}
The diagonal sum $\diag(\df,n)$ equals the number of partitions of $n$ with no part 2.
\end{cor}

\begin{proof}
As in Corollary \ref{excdiag}, letting $w=q$ in \eqref{dfc2} gives the diagonal sums of the $\df(n,k)$ triangle.  The resulting ratio
\[ \frac{(1/q;q)_n}{(q;q)_n}\]
in the summation is nonzero only for $n = 0$ and 1, leading to 
\[\frac{1+q}{(q^2;q)_\infty} = \frac{1-q^2}{(q;q)_\infty}\]
where the latter is a generating function for partitions with no part 2.
\end{proof}

We will focus on more combinatorial aspects of partition deficiencies in future work.

\end{document}